\documentclass[leqno]{amsart}

\usepackage{latexsym,amssymb,amsthm,amsmath,amscd,xcolor}

\theoremstyle{plain}
\newtheorem{theorem}{Theorem}[section]

\numberwithin{equation}{section}

\theoremstyle{remark}
\newtheorem{remark}{Remark}[section]

 \numberwithin{equation}{section}

\newtheorem*{Theorem A}{{\bf Theorem A}}
\newtheorem*{Theorem B}{{\bf Theorem B}}
\newtheorem*{Theorem C}{Theorem C}

 \numberwithin{equation}{section}
\def\<{\left < }
\def\>{\right >}
\def\({\left ( }
\def\){\right )}
\def\e{\eqref}

\begin{document}

\title[Some links between $F$-harmonicity, submersion and cohomology] {Some links between $F$-harmonicity, submersion and cohomology}

\author[B.-Y. Chen]{Bang-Yen Chen}
\address{Department of Mathematics\\
    Michigan State University \\East Lansing, Michigan 48824--1027\\ U.S.A.}
\email{chenb@msu.edu}

\author[S. W. Wei]{Shihshu Walter Wei}
\address{Department of Mathematics\\
University of Oklahoma\\ Norman, Oklahoma 73019-0315\\ U.S.A.}
\email{wwei@ou.edu}

\begin{abstract}  By studying cohomology classes that are related with $p$-harmonic morphisms, $F$-harmonic maps, and $f$-harmonic maps, we extend several of our previous results on Riemannian submersions and $p$-harmonic morphisms to $F$-harmonic maps, and $f$-harmonic maps which are submersions.
\end{abstract}


 
\thanks{}

\date{}

\maketitle

\section{Introduction}\label{S1}

We note in the study of $p$-harmonic maps between Riemannian manifolds that the following duality occurs between  (A) conformal immersions and (B) horizontally conformal submersions: 
\vskip.1in

\noindent 
  { (A) {\bf Theorem} (\cite {T})   {\it Let $u: (M,g_{M}) \to (N,g_{N})$ be a conformal immersion with $\dim M=p$. Then 
{\it $u$ is a $p$-harmonic map  if and only if $u(M)$ is a (parametric) minimal submanifold in $N$.}}
\vskip.1in

\noindent 
(B)  { {\bf Theorem} (\cite {BG}) {\it Let $u: (M,g_{M}) \to (N,g_{N})$ be a horizontally conformal submersion   with $\dim N=p$. Then  $u$  is a $p$-harmonic map if and only if all the fibers of $u$ are (parametric) minimal submanifolds in $M$.}}}

\vskip.1in
\noindent  This duality between (A) and (B) links $p$-harmonic maps and (parametric) minimal submanifolds.  

On the other hand, there also exists a duality between (C) (non-parametric) minimal hypersurface in the Euclidean space $\mathbb R^{n+1}$ and (D) maximal spacelike (non-parametric) hypersurface in  the Minkowski space $\mathbb R^{n,1}$ equipped with the Lerentzian metric $ds^2 =\sum_{i=1}^n (dx^i)^2 -dt^2 $, where 

\vskip.1in

\noindent (C) { The} minimal hypersurface in the Euclidean space $\mathbb R^{n+1}$  is given by the graph of a function $u$ on a Euclidean domain satisfying the differential equation
\begin{equation}\operatorname{div}\bigg ( \frac {\nabla u}{\sqrt{1+ |\nabla u|^2}} \bigg )=0; \label{1.1}\end{equation}

\vskip.1in
and 

\vskip.1in
\noindent (D) { The} maximal spacelike hypersurface in the Minkowski space $\mathbb R^{n,1}$ is furnished by the graph of a function $v$ on a Euclidean domain satisfying a dual differential equation
\begin{equation}\operatorname{div}\bigg ( \frac {\nabla v}{\sqrt{1- |\nabla v|^2}} \bigg )=0,\;\; \;\;  |\nabla v|^2<1. \label{1.2}\end{equation}

 E. Calabi showed in \cite{C} that equations \eqref{1.1} and \eqref{1.2} are equivalent over any simply-connected domain in $\mathbb R^2$.
Furthermore, the notion of $F$-harmonic maps { unifies}  minimal hypersurfaces and  maximal spacelike hypersurfaces;  namely, the solutions $u$ and $v$ are $F$-harmonic maps from a domain in {  $\mathbb R^{n+1}$} to $\mathbb R$ with $F(t) = \sqrt{1+2t} -1$ and $F(t) = 1-\sqrt{1-2t}$, respectively (see \cite {W9,DW}). Hence, one can unify harmonic maps, $p$-harmonic maps, $\alpha$-harmonic maps, exponential harmonic maps, (nonparametric) minimal hypersurfaces in Euclidean space, and (nonparametric) maximal spacelike hypersurface in Minkowski space $\mathbb R^{n,1}$ by $F$-harmonic maps (see \cite {A, W9}). Also, the ``$F\, $"-idea can be extended to Gauge Theory (see \cite {DW}). Here
\begin{equation} F : [0,\infty) \to [0,\infty)\quad  \operatorname{with}\quad F(0) = 0\quad \operatorname{is}\, C^2\,  \operatorname{strictly}\,  \operatorname{increasing}. \label{1.3}\end{equation} 

We recall in the study of topology on a compact Riemannian manifold $M$, it is well-known that nontrivial fundamental groups $\pi_1(M)\, ,$ homology groups, and cohomology classes can be represented by stable closed geodesics,  stable minimal rectifiable currents, and harmonic forms on $M$,  respectively, by Cartan's Theorem (\cite {Ca}),  Federer-Fleming's Theorem (\cite {FF}), and Hodge Theorem (\cite {H}).  In an analogous spirit, it was shown in \cite {S.W.Wei 1} that homotopy classes can be represented by $p$-harmonic maps
 (For definition and examples of $p$-harmonic maps, see e.g. \cite {W9}):
\vskip.1in

\noindent {\bf Theorem A.} {\it If $N$ is a compact Riemannian manifold, then for any positive integer $i$, each class in the i-th homotopy group $\pi_i(N)$ can be represented by a $C^{1,\alpha}$ $p$-harmonic map $u_0$ from an $i$-dimensional sphere $S^i$ into $N$ minimizing $p$-energy in its homotopy class for any $p>i$.}
\vskip.1in

Further applications and homotopically vanishing theorems were
explored in \cite {S.W.Wei 1,WW}, using stable $p$-harmonic maps as catalysts, whereas homologically vanishing theorem was
given in \cite {LS, W2, HoW}, using stable rectifiable currents as catalysts. 

On the other hand, B.-Y. Chen established in \cite{c2005} the following result involving Riemannian submersion, minimal immersion, and  cohomology class.
\vskip.1in

\noindent {\bf Theorem B.} (\cite{c2005}) {\it Let $b = \dim \, B$ and let $\pi : (M,g_{M})\to
(B,g_{B})$ be a Riemannian submersion with minimal fibers and orientable
base manifold $B$. If  $M$ is a closed manifold with cohomology
class $H^{b}(M,\mathbf R) =0$, then the horizontal distribution
$\mathcal{H}$ of the Riemannian submersion is never integrable.
Thus the submersion $\pi$ is never non-trivial.}
\vskip.1in

Whereas $p$-harmonic maps represent homotopy classes,  B.-Y. Chen and S. W. Wei connected  the two seemingly unrelated areas of $p$-harmonic morphisms and cohomology classes in the following.
\vskip.1in

\noindent {\bf Theorem C.} {\color{red} (\cite {CW09,CW10})} {\it Let $u: (M,g_{M}) \to (N,g_{N})$ be an $n$-harmonic morphism which is a submersion with $\dim N=n$. If $N$ is orientable and $M$ is a closed manifold with $n$-th cohomology class $H^n(M,\mathbf R) = 0$, then the horizontal distribution $\mathcal H$ of $u$ is never integrable. Hence the submersion $u$ is always non-trivial.}
\vskip.1in

This recaptures Theorem B when $\pi : M\to B$ is  a  Riemannian submersion with minimal fibers
and orientable base manifold $B$. While a horizontally weak conformal $p$-harmonic map is a
$p$-harmonic morphism (see Theorem \ref{T:2.5}), $p$-harmonic morphism is also linked to cohomology class as follows.
\vskip.1in
\noindent {\bf Theorem D.} {\color{red} (\cite {CW09,CW10})} {\it Let $u: (M,g_{M}) \to (N,g_{N})$ be an $n$-harmonic morphism  with $ n=\dim  N$ which is a submersion. Then the pull back of the volume element of $N$ is a harmonic $n$-form if and only if the horizontal distribution $\mathcal H$ of $u$ is completely integrable.}
\vskip.1in

Following the proofs given in \cite {CW09,CW10}, and by applying a characterization theorem of a $p$-harmonic morphism (see Theorem 2.5) and Theorem 2.6, we obtain a dual version of Theorem D. In particular, we have $p$-harmonic maps and cohomology classes are interrelated in \cite {W9} as follows.
\vskip.1in

\noindent {\bf Theorem E.} {\it  Let $M$ be a closed manifold, and $u : M \to N$ be an $n$-harmonic map with $n = \dim N$ which is a submersion. Assume that the horizontal distribution $\mathcal H$ of $u$ is integrable and $u$ is an $n$-harmonic morphism. Then we have $H^n(M,\mathbf R) \ne 0$.}
\vskip.1in

\noindent {\bf Theorem F.} (\cite {W9}) {\it Let $u : M \to N$ be an $n$-harmonic map with $n = \dim N$ which is a submersion. Let the horizontal distribution $\mathcal H$ of $u$ be integrable. If $M$ is a closed manifold with cohomology class $H^n(M,\mathbf R) = 0$. Then $u$ is not an $n$-harmonic morphism. Thus the submersion $u$ is always nontrivial.}

\section{Preliminaries}\label{S2}

We recall some basic facts, notations,  definitions, and formulas
for minimal submanifolds, submersions, $p$-harmonic morphisms, $F$-harmonic maps, and $f$-harmonic maps for later use (see, e.g.,  {\color{red}\cite{book73,c2005,book11,CW09}} for details).

\subsection{Basic formulas and equations}\label{S2.1}

Let $\widetilde M$ be a Riemannian manifold with Levi-Civita
connection $\widetilde \nabla$.  The tangent bundle of $\widetilde M$ is
denoted by $T \widetilde M$, and the (infinite dimensional) vector
space of smooth sections of a smooth vector bundle $E$ is denoted
by $\Gamma (E)$. Let $M$ be a submanifold of dimension $n\geq 2$
in $\widetilde M$. Denote by $\nabla$ and $D$, the Levi-Civita
connection and the normal connection of $M$, respectively. For each
normal vector $\xi\in T^\perp_xM,\, x\in M$, the shape operator
$A_{\xi}$ is a symmetric endomorphism of the tangent space $T_xM$
at $x$. Then the shape operator and the second fundamental form
$\mathsf h$ are related by
\begin{align}\label{2.1} \<\mathsf h(X,Y),\xi\>=\<A_{\xi}X,Y\>\end{align}
for $X,Y$ tangent to $M$ and $\xi$ normal to $M$.   

The formulas of Gauss and Weingarten are given respectively by (cf.
\cite{book73})
\begin{align} &\label{2.2}\widetilde \nabla_XY=\nabla_XY+ \mathsf h(X,Y), 
\\&\label{2.3}\widetilde \nabla_X \xi=-A_\xi X+D_X\xi\end{align} for tangent vector
fields $X,Y$ and normal vector field $\xi$ on $M$. 

 The \emph{mean curvature vector field} of a submanifold $M$
is defined by $H=\frac{1}{n}{\rm trace}\, \mathsf h$. A submanifold $M$ in
$\widetilde M$ is called \emph{totally geodesic} (respectively,
\emph{minimal}) if its second fundamental form $\mathsf h$ (respectively,
its mean curvature vector field $H$)  vanishes identically.

\subsection{Energy and harmonic maps}\label{S2.2}
Let $u: (M,g_{M}) \to (N,g_{N})$ be a differential map between two Riemannian manifolds $M$ and $N$. Denote  $e_u$ the {\it energy density} of  $u$,  which is given by 
\begin{equation}\label{2.4} e_u = \frac 12\sum _{i=1}^m g_{N}\big (du(e_i),du(e_i)\big ) = \frac 12 |du|^2\, , \end{equation}
where $\{e_1, \cdots, e_n\}$ is a local orthonormal frame field on $M$ and $|du|$ is the Hilbert-Schmidt norm of $du$, determined by the metric $g_{M}$ of $M$ and the metric $g_{N}$ of $N$. 
{\it The energy of $u$}, denoted by $E(u)$,  is defined to be $$E(u) = \int_M e_u\, dv_g.$$ A smooth map $u: M \to N$ is called {\it harmonic} if $u$ is a critical point of the energy functional $E$
with respect to any compactly supported variation.

\subsection{Submersions}\label{S2.3}

A differential map $u: (M,g_{M}) \to (N,g_{N})$  between two Riemannian manifolds  is called a submersion at a point $x\in M$ if its differential $du_{x}: T_x(M)\to T_{u(x)}(N)$ is a surjective linear map. A differentiable map $u$ that is a submersion at each point $x\in M$ is called a {\it submersion}. For each point $x\in N$, $u^{-1}(x)$ is called a {\it fiber}.
For a submersion $u: M \to N$, let $\mathcal H_x$ denote the orthogonal complement of Kernal$\, \big (du_x : T_x(M)\to T_{u(x)}(N)\big )$ in $T_x(M)$ at $x\in M$. Let $\mathcal H=\{\mathcal H_x: x\in M\}$ denote the horizontal distribution of $u$. 
A submersion $u: (M,g_{M}) \to (N,g_{N})$ is called {\it horizontally weakly conformal} if the restriction of $du_x$ to $\mathcal H _x$ is conformal, i.e., there exists a smooth function $\lambda$ on $M$ such that 
\begin{align}\label{2.5} u^{*} g_{N} = \lambda^2  {g _M}_ {|_{\mathcal H}} \;\;\;\; {\rm or}\;\; \;\;
g_{N}\big (du_x(X),du_x(Y)\big) = \lambda ^2 (x) g_{M}(X, Y)\end{align}
 for all $X,Y \in \mathcal{H} _x$ and $x \in M$. 
  If the function $\lambda$ in \eqref{2.5} is  positive, then $u$ is called  {\it horizontally conformal}
and   $\lambda$ is called the {\it dilation} of $u$.  For a horizontally conformal submersion $u$ with dilation $\lambda$, the {\it energy density} is  $e_u = \frac{1}{2} {n\lambda ^2}.$ 
A horizontally conformal submersion with dilation $\lambda \equiv 1$ is called a {\it Riemannian submersion}. 
For Riemannian submersions, we have the following results.

\begin{theorem}[\cite{c2005}, Theorem 2]\label{T:2.1} Let $\pi : (M,g_{M})\to (B,g_{B})$ be a Riemannian submersion from a closed manifold $M$ onto an orientable base manifold $B$. Then the pullback of the volume element of $B$ is harmonic if and only if the horizontal distribution $\mathcal{H}$ is integrable and  fibers are minimal.
\end{theorem}

\begin{theorem}[\cite {W9}, Theorem 2.5] \label{T:2.2} Let $\pi : (M,g_{M})\to (B,g_{B})$ be a Riemannian submersion. Then $\pi$ is a $p$-harmonic map for every $p > 1$ if and only if all fibers $\pi^{-1}(y) , y \in B$, are minimal submanifolds in $M$.\end{theorem}

\subsection{Harmonic morphisms}\label{S2.4}

A $C^2$--map $u: (M,g_{M}) \to (N,g_{N})$ is said to be a \emph{harmonic morphism} if for any harmonic function $f$ defined on an open set $V$ of $N$, the composition $f \circ u$ is harmonic on $u^{-1}(V)$. 

P. Baird and J. Eells used the stress-energy tensor to establish the following.

\begin{theorem}[\cite{BE}] \label{T:2.3} Let $u: (M,g_{M}) \to (N,g_{N})$ be a harmonic morphism which is a submersion everywhere on $M$. Then $\big ($setting $n = \dim N$, and $e_u$ as in \eqref{2.4}$\big )$

\begin{itemize}
\item[{\rm (a)}]  if $n= 2$, the fibres are minimal submanifolds ;

\item[{\rm (b)}] if $n  >  2$ , then the following properties are equivalent :
\begin{itemize}
\item[{\rm (i)}] the fibres are minimal submanifolds;
 
\item[{\rm (ii)}] $\operatorname{grad} (e_u)$ is a vertical field;

\item[{\rm (iii)}] the horizontal distribution has mean curvature vector $\frac {\operatorname{grad}(e_u)}{2 e_u}$.
\end{itemize}\end{itemize}
\end{theorem}

\subsection{$p$-harmonic morphisms}\label{S2.5}

A $C^2$--map $u: (M,g_{M}) \to (N,g_{N})$ is said to be a \emph{$p$-harmonic
morphism} if for any $p$-harmonic function $f$ defined on an open
set $V$ of $N$, the composition $f \circ u$ is $p$-harmonic on
$u^{-1}(V)$. 
We have the following link between $p$-harmonic morphisms for every $p > 1$ and minimal fibers:

\begin{theorem}[\cite {W9}, Proposition 2.4] \label{T:2.4} If $\pi: (M,g_{M}) \to (N,g_{N})$ is a Riemannian submersion, then $\pi$ is a $p$-harmonic morphism for every $p > 1$ if and only if all fibers $\pi^{-1}(y),\, y \in N$, are minimal submanifolds of $M$. \end{theorem}

E. Loubeau \cite{Lo} and J. M. Burel and E. Loubeau \cite{BL} obtained the following characterization of $p$-harmonic morphisms.

\begin{theorem}[\cite{Lo, BL}] \label{T:2.5}  A $C^2$--map $u: (M,g_{M}) \to (N,g_{N})$ is a $p$-harmonic morphism with $p \in (1, \infty)$ if and only if it is a $p$-harmonic, horizontally conformal map. \end{theorem}

Theorem 2.4 (\cite {W9}, Proposition 2.4) 
links $p$-harmonic maps for every $p >1$ with minimal fibers in the presence of Riemannian submersion. This is dual to a minimal submanifold occurs from a $p$-harmonic map for every $p >1$ in the presence of isometric immersion 
(cf. \cite {S.W.Wei 1} ).  In \cite {BG}, P. Baird and S. Gudmundsson linked $n$-harmonic morphisms with minimal fibers as follows.

\begin{theorem} $($\cite {BG}$)$\label{T:2.6}  If $u: (M,g_{M}) \to (N,g_{N})$ is a horizontally conformal submersion from a manifold $M$ onto a manifold $N$ with $n=\dim N$, then $u$ is $n$-harmonic if and only if the fibers of $u$ are minimal in $M$. \end{theorem}

\subsection{$F$-harmonic maps}\label{S2.6}

Let $F$ be as in \eqref{1.3} and let $u: (M,g_{M}) \to (N,g_{N})$ be a smooth map between two compact Riemannian manifolds. Then the map $u:M\to N$ is called {\it $F$-harmonic} if it is a critical point of the $F$-energy functional: 
$E_F (u) =\int _M F\! \(\frac {|du|^2}{2}\) dv_g.$  In particular,
when
\begin{equation*}\begin{aligned}\label{2.7}F (t) = &\, \; t,\; \; \frac 1p (2t)^{\frac p2}, \;\; (1 + 2t)^{\alpha} - 1\, \;\; (\alpha > 1, \dim M =2),
 \\& e^t - 1,\;\;  \sqrt {1 + 2t} - 1, \;\;  1- \sqrt {1 - 2t},\;\;  {\rm respectively,} 
 \end{aligned}\end{equation*}
the $F$-energy $E_F(u)$ becomes energy, $p$-energy, (normalized) $\alpha$-energy,  (normalized) exponential energy, (normalized) area functional in Euclidean space, and (normalized) area functional in Minkowski space. Hence, its critical point $u$ or its graph is harmonic, $p$-harmonic, $\alpha$-harmonic, exponential harmonic, minimal hypersurface in Euclidean space $\mathbb R^{n+1}$, maximal spacelike hypersurface in Minkowski space $\mathbb R^{n+1}_1$, respectively. 
 
 For horizontally conformal $F$-harmonic maps, M. Ara proved the following.

\begin{theorem}  [cf. \cite {A}] \label{T:2.7}  Let $u: (M,g_{M}) \to (N,g_{N})$, $m > n$, be an $F$-harmonic map, which is horizontally conformal  with dilation $\lambda$. Assume that the zeros of $(n-2) F^{\prime}(t)-2tF^{\prime\prime}(t)$ are isolated. Then the following three properties are equivalent:
\begin{itemize} 
\item[{\rm (1)}] The fibers of $u$ are minimal submanifolds.
 \item[{\rm (2)}]  $\operatorname{grad}(\lambda^2)$ is vertical.
\item[{\rm (3)}] The horizontally distribution of $u$ has mean curvature vector $\frac{\operatorname{grad}(\lambda^2)}{2\lambda^2}$.
\end{itemize}
\end{theorem}

\subsection{$f$-harmonic maps}\label{S2.7}

Let $f : (M,g) \to (0,\infty)$ be a smooth function. The notion of $f$-harmonic maps was first introduced and studied by A. Lichnerowicz in 1970 (see \cite{L}). 
A map $u: (M,g_{M}) \to (N,g_{N})$ between two Riemannian manifolds is said to be {\it $f$-harmonic}, if $u$ is a critical point of the {\it $f$-energy functional} $E_f$ with respect to any compactly supported variation (cf. \cite {EL, L}), where
\begin{equation}E_f(u)= \frac 12 \int _M f|du|^2\, dv_g\, .\end{equation}
Examples of $f$-harmonic maps include harmonic maps with $f \equiv C\, ,$ a positive constant, and {\it submersive $p$-harmonic maps} with $f=|du|^{p-2}$. 

The next result is due to Y.-L. Ou.

\begin{theorem}[\cite{Ou}] \label{T:2.8} Let $u: (M,g_{M}) \to (N,g_{N})$ be a horizontally weakly conformal map with dilation $\lambda$. Then, any two of the following conditions imply the other one:
\begin{itemize}
\item[{\rm (a)}]  $u$ is an $f$-harmonic map;
\item[{\rm (b)}]  $\operatorname{grad}( f \lambda^{2-n})$ is vertical;
\item[{\rm (c)}]  $u$ has minimal fibers. \label{T:4.1}
\end{itemize}\end{theorem}

\section{Statements Of Theorems}\label{S3}

In the following we shall assume that $\dim M=m> n=\dim N$.
The purpose of this article is to extend the ideas in \cite {CW09, W9} and to connect the seemingly unrelated areas of $F$-harmonic maps, $f$-harmonic maps, and cohomology classes. 
More precisely we prove the following
\vskip.1in

\noindent{\bf Theorem 1.}  {\it Let $u: M \to N$ be an $F$-harmonic map which is a horizontally conformal submersion with dilation $\lambda$. And let $\operatorname{grad}$ $\lambda^2$ be vertical, and  the zeros of $(n-2) F^{\prime}(t)-2tF^{\prime\prime}(t)$ be isolated. If $H^n(M,\mathbf R) = 0$, then the horizontal distribution $\mathcal{H}$ of $u$ is not integrable. Thus the submersion is always non-trivial.}

\begin{remark}\label{R:3.1}  $(i)\, $Theorem 1 generalizes Theorem B, in which by assumption $u$ is a horizontal conformal submersion with dilation $\lambda =1$, and $u$ is a $p$-harmonic map for every $p > 1$ according to  Theorem \ref{T:2.4}. If we choose  $p \ne n$, then $\operatorname{grad}$ $\lambda^2 = 0\, ,$
$F(t) = \frac 1p (2t)^{\frac p2}$ so that  $F$-harmonic map becomes $p$-harmonic map, $p \ne n$, and the zeros of $(n-2) F^{\prime}(t)-2tF^{\prime\prime}(t)$ are 
isolated. 

\vskip.05in
\noindent
$(ii)$ Theorem 1 augments Theorem C but does not generalize Theorem C, since when $F(t) = \frac 1n (2t)^{\frac n2}$, $F$-harmonic map $u$ becomes $n$-harmonic map, and the zeros of $(n-2) F^{\prime}(t)-2tF^{\prime\prime}(t)$
are not isolated.  Thus $n$-harmonic maps are not applicable to Theorem 1. 

\vskip.05in
\noindent
$(iii)$ Theorem 1 is equivalent to the following Theorem 2, due to the equivalence in Theorem \ref{T:2.7}. 
\end{remark}		

\noindent{\bf Theorem 2.}  {\it Let $u: M \to N$ be an $F$-harmonic map which is a horizontally conformal submersion with dilation $\lambda$. Assume that the horizontal distribution
has mean curvature vector $\frac{\operatorname{grad}(\lambda^2)}{2\lambda^2}$, and  the zeros of $(n-2) F^{\prime}(t)-2tF^{\prime\prime}(t)$ are isolated. If $H^n(M,\mathbf R) = 0$, then the horizontal distribution $\mathcal{H}$ of $u$ is not integrable. Thus the submersion is always non-trivial.}
\vskip.05in

\noindent
\noindent{\bf Theorem 3.}  {\it Let $u: M \to N$ be an $F$-harmonic map which is a horizontally conformal submersion with dilation $\lambda$. Assume that the horizontal distribution
has mean curvature vector $\frac{\operatorname{grad}(\lambda^2)}{2\lambda^2}$, and  the zeros of $(n-2) F^{\prime}(t)-2tF^{\prime\prime}(t)$ are isolated. Then  the pullback of the volume element of $N$ is a harmonic $n$-form if and only if the horizontal distribution  $\mathcal{H}$ of $u$ is completely integrable.}

\begin{remark}\label{R:3.2}  By similar arguements as in Remark \ref{R:3.1},  $(i)\,$ Theorem 2  generalizes Theorem 2.1;
$(ii)$ Theorem 2  augments Theorem D but does not generalize Theorem D; and 
$(iii)$ Theorem 2  is equivalent to the following Theorem 4. 
\end{remark}		

\noindent{\bf Theorem 4.}  {\it Let $u: M \to N$ be an $F$-harmonic map which is a horizontally conformal submersion with dilation $\lambda$. Assume that $\operatorname{grad}$ $\lambda^2$ be vertical, and  the zeros of $(n-2) F^{\prime}(t)-2tF^{\prime\prime}(t)$ are isolated. Then  the pullback of the volume element of $N$ is a harmonic $n$-form if and only if the horizontal distribution  $\mathcal{H}$ of $u$ is completely integrable.}
\vskip.05in

We then study  $f$-harmonic maps which was first introduced and studied by A. Lichnerowicz in 1970 $($see \cite{L}$)$, where $f : (M,g) \to (0,\infty)$ is a smooth function.. 
Examples of $f$-harmonic maps include harmonic maps with $f \equiv C\, ,$ a positive constant, and {\it submersive $p$-harmonic maps} with $f=|du|^{p-2}$ (cf. \S4).
We have the following results.
\vskip.05in

\noindent{\bf Theorem 5.} {\it Let  $u: M \to N$ be an $f$-harmonic map which is a horizontally conformal submersion with dilation $\lambda$  such that $\operatorname{grad} (f \lambda^{2-n})$ is vertical. If $M$ is a closed manifold with cohomology class $H^n(M,\mathbf R) = 0$, then the horizontal distribution $\mathcal H$ of $u$ is never integrable. Thus the submersion $u$ is always non-trivial.}
\vskip.05in

This recaptures Theorem B, when $f$-harmonic map is an $n$-harmonic map, then $f = |du|^{n-2}$.
\vskip.05in

\noindent{\bf Theorem 6.} {\it Let  $u: M \to N$ be an $f$-harmonic map which is a horizontally conformal submersion and with dilation $\lambda$  such that $\operatorname{grad} (f \lambda^{2-n})$ is vertical. Then the pullback of the volume element of $N$ is a harmonic $n$-form if and only if  horizontal distribution is completely integrable.}
\vskip.05in

This recaptures Theorems D, when $f$-harmonic map is an $n$-harmonic map.

\section{Proofs of Theorems 1-6}

\begin{proof}  {\bf Proof of Theorem 3:} Let $u: (M,g_{M}) \to (N,g_{N})$ be an $F$-harmonic map which is a horizontally conformal submersion with dilation $\lambda$ and $n=\dim N$. Assume that the horizontal distribution
has mean curvature vector $\frac{\operatorname{grad}(\lambda^2)}{2\lambda^2}$, and  the zeros of $(n-2) F^{\prime}(t)-2tF^{\prime\prime}(t)$ are isolated. Then
the submersion $u$ has minimal fibers according to Theorem \ref{T:2.7}.

Let $\{{\bar e}_{1},\ldots,{\bar e}_{n}\}$ be an oriented local orthonormal frame of the base manifold $(N,g_{N})$ and let ${\bar \omega}^{1},\ldots,{\bar \omega}^{n}$ denote the dual 1-forms of $\{{\bar e}_{1},\ldots,{\bar e}_{n}\}$ on $N$. Then $${\bar \omega}={\bar \omega}^{1}\wedge \cdots  {\bar \omega}^{n}$$ is the volume form of $(N,g_{N})$, which is a closed $n$-form on $N$. 

Consider the pull back of the volume form ${\bar \omega}$ of $N$ via $u$, which is denoted by  $u^{*}({\bar \omega})$.  Then $u^{*}({\bar \omega})$ is a simple $n$-form on $M$  satisfying
  \begin{equation}\label{4.1}d \big (u^{*}({\bar \omega})\big )=u^{*}(d{\bar \omega}) =0,\end{equation} since the exterior differentiation $d$ and the pullback $u^{*}$ commute. 
  
Assume that $\dim M =m= n+k$ and let $e_1, \dots, e_{n+k}$ be a local orthonormal frame field with $\omega^1, \dots, \omega^{n+k}$ being its dual coframe fields  on $M$ such that 
  
\vskip.06in
\noindent   (i) $e_{1}, \dots, e_{n}$ are basic horizontal vector fields satisfying $du(e_{i})=\lambda \bar e_{i}$, $ i=1,\ldots,n$, and  $du(e_{1}), \dots, du(e_{n})$ give a positive orientation of $N$; and
  
 \noindent  (ii) $e_{n+1}, \dots, e_{n+k}$ are vertical vector fields. 

\vskip.06in
\noindent Then  we have 
\begin{equation}\label{4.2}  \omega^j(e_s) = 0,\;\;\;  \omega^i(e_j)=\delta_{ij},\;\; 1 \le i, j \le n; \; n+1 \le s \le n+k\, . \end{equation} 
Also, it follows from (i) that
\begin{equation}\label{4.3}  u^{\ast} \bar{\omega}^{i}=\frac{1}{\lambda}\omega^{i},\;\; \; i=1,\ldots,n.
 \end{equation} 
  
 If we put 
\begin{equation}\label{4.4} \omega = \omega^1 \wedge \cdots \wedge \omega^{n}\;\;\; {\rm and}\;\;\; \omega^{\bot} = \omega^{n+1}\wedge \cdots \wedge \omega^{n+k}\, , \end{equation}
then 
\begin{equation}\label{4.5}d\omega^{\bot} = \sum _{i=1}^{k} (-1)^{i}\omega^{n+1}\wedge \dots \wedge d\omega^{n+i}\wedge  \cdots \wedge \omega^{n+k}.
\end{equation}
It follows from \e{4.2} and \e{4.5}  that $d\omega^{\bot}=0$ holds identically if and only if the following two conditions are satisfied:
\begin{equation} \label{4.6}d\omega^{\bot}(e_i,e_{n+1}, \dots, e_{n+k}) = 0 \end{equation}
for $i=1,\ldots, n$; and
\begin{equation} \label{4.7}d\omega^{\bot}(X,Y,V_1,\dots,V_{k-1}) = 0. \end{equation}
for any horizontal vector  fields $X,Y$ and vertical vector fields $V_{1}, \ldots, V_{k-1}$.

Since the fibers of $u$ are minimal in $M$,  for each $1 \le i \le n$, we find
\begin{equation}\begin{aligned}\label{4.8}
d\omega^{\bot}(e_i,\,&e_{n+1}, \dots, e_{n+k})\\
& = \sum _{j=1}^k (-1)^{j+1} \omega^{\bot}([e_i,e_{n+j}], e_{n+1},\dots,\hat{e} _{n+j}, \dots, e_{n+k}) \\
&=\sum _{j=1}^k (-1)^{j+1} \big(\omega^{n+j} (\nabla _{e_i} e_{n+j}) - \omega^{n+j} (\nabla _{e_n+j} e_{i}) \big)\\
&= \sum _{j=1}^k - \langle \nabla _{e_{n+j}} e_{i}, e_{n+j}\rangle\\
&=  \sum _{j=1}^k \langle h(e_{n+j}, e_{n+j}), e_i\rangle \\
&=  0, \end{aligned}
\end{equation}
where ``$\,\hat{\cdot}\,$'' denotes the missing term, which proves that condition \e{4.6} holds.

Next, we suppose that the horizontal distribution $\mathcal{H}$ is integrable. If $X, Y$ are  horizontal vector fields,   then $[X,Y]$ is also horizontal by Frobenius theorem. So, for vertical vector fields $V_1, \dots, V_k\, ,$ we find (cf. \cite[formula (6.7)]{c2005} or \cite[formula (3.5)] {W9})
\begin{equation} \label{4.9}d\omega^{\bot}(X,Y,V_1,\dots,V_{k-1}) = \omega^{\bot}([X,Y],V_1,\dots,V_{k-1}) = 0. \end{equation}
Consequently, from \e{4.8} and \e{4.9} we get
\begin{equation}\label{4.10}d\omega^{\bot}=0.\end{equation}  

Next, we  show that if  $\mathcal{H}$ is integrable, then 
we have $d (u^{\ast} \bar{\omega}\,)^{\bot}=0\, .$ 
Since $u$ is a horizontally conformal submersion with dilation $\lambda$, it preserves orthogonality, which is crucial to horizontal and vertical distributions, and the pullback $u^{\ast}$ expands the length of $1$-form by $\frac 1\lambda$ in every direction. This, via \eqref{4.3} and \eqref{4.10},  leads to
\begin{equation}\begin{aligned}\label{4.11}
d(u^{\ast} \bar{\omega}\,)^{\bot}& = d(u^{\ast} \bar{\omega}^{1} \wedge \dots \wedge u^{\ast} \bar{\omega}^{n})^{\bot}  \\
& = d\left(\frac 1{\lambda} {\omega}^{1} \wedge \dots \wedge \frac 1{\lambda} {\omega}^{n}\right)^{\bot}  \\
& = d\left(\frac 1{\lambda} \omega^{n+1} \wedge \dots \wedge \frac 1{\lambda}  \omega^{n+k}\right)  \\
& = d\left(\frac 1{\lambda^k} \omega^{\bot}\right)  \\
&=  -\frac{k}{2} \lambda ^{-k-2} d (\lambda ^2) \wedge \omega^{\bot} + \frac 1{\lambda^k} d \omega^{\bot}
\\ &=  -\frac{k}{2} \lambda ^{-k-2} d (\lambda ^2) \wedge \omega^{\bot}
\\& =0.\end{aligned}\end{equation}
The last equality  of \eqref{4.11} is due to the fact that grad$(\lambda^2)$ is vertical and Theorem 2.7. Therefore,  grad$(\lambda^2)$, $e_{n+1}, \cdots, e_{n+k}\, $ are linearly dependent. Hence, we obtain $d (\lambda ^2) \wedge \omega^{\bot}=0$.
Because $d(u^{\ast}\omega)^{\bot}=0$ is equivalent to $u^{\ast}\omega$ being co-closed, it follows from $d(u^{\ast}\omega)=0$ and $d(u^{\ast}\omega)^{\bot}=0$ that the pullback of the volume form, $u^{\ast}\omega$, is a harmonic $n$-form on $M$ whenever $\mathcal H$ is completely integrable.

Conversely, it follows from the proof given above  that if $u^{\ast}\omega$ is a harmonic $n$-form on $M$, then the horizontal distribution $\mathcal{H}$ is completely integrable.
This proves Theorem 3.
\vskip.05in
  
\noindent
{\bf Proof of Theorem 1:} Since each nonzero harmonic form represents a nontrivial cohomology class by Hodge theorem (\cite {H}), and since on a closed manifold a differential form is harmonic if and only if it is closed and co-closed, Theorem 1 follows from Theorem 3.
\vskip.05in

\noindent
{\bf Proof of Theorem 2 and Theorem 4:} Theorem 2 and Theorem 4 follow from Theorem 1 together with Remarks \ref{R:3.1} and \ref{R:3.2} respectively.
\vskip.05in

\noindent
{\bf Proof of Theorem 5 and Theorem 6:} By virtue of Theorem \ref{T:2.8}, the submersion $u$ in Theorem 5 and Theorem 6 has minimal fibers, we proceed as before, $u$ is an $n$-harmonic map by virtue of Theorem \ref{T:2.6}, and hence $u$ is an $n$-harmonic morphism according to Theorem \ref{T:2.5}.  We are then ready to apply Theorem C and Theorem D to complete the proof of Theorem 5 and Theorem 6 , respectively.
\end{proof}

\end{document}